\documentclass{amsart}

\usepackage{amssymb,latexsym}

\numberwithin{equation}{section}

\newtheorem{theorem}{Theorem}[section]
\newtheorem{proposition}[theorem]{Proposition}

\newtheorem{corollary}[theorem]{Corollary}
\newtheorem{remark}[theorem]{Remark}
\newtheorem{lemma}[theorem]{Lemma}
\newtheorem{example}[theorem]{Example}
\newtheorem{definition}[theorem]{Definition}

\def\proof{\smallskip\noindent {\bf Proof. }}
\def\endproof{\hfill$\square$\medskip}

\def\dimv{\underline{\dim}}

\def\ZZ{\mathbb{Z}}

\def\RR{\mathbb{R}}

\newcommand{\ses}[3]
{\mbox{$0 \rightarrow #1 \rightarrow #2 \rightarrow #3 \rightarrow 0$}}

\newcommand{\drc}{\rep\widetilde{\Gamma}}
\newcommand{\tdimv}{\underline{\textrm{sdim}}}
\newcommand{\codim}{\mbox{codim }}
\newcommand{\rep}{{\rm rep }\,}
\newcommand{\vect}{{\rm vect }}

\begin{document}

\title[Generalized associahedra via quiver representations]
{Generalized associahedra via quiver representations}

\author{Robert Marsh}
\address{Department of Mathematics and Computer Science, University of
Leicester, University Road, Leicester LE1 7RH, England}
\email{R.Marsh@mcs.le.ac.uk}

\author{Markus Reineke}
\address{BUGH Wuppertal, Gau\ss stra\ss e 20, D-42097 Wuppertal, Germany}
\email{reineke@math.uni-wuppertal.de}

\author{Andrei Zelevinsky}
\address{Department of Mathematics, Northeastern University,
  Boston, MA 02115, USA}
\email{andrei@neu.edu}

\begin{abstract}
We provide a quiver-theoretic interpretation of certain smooth complete
simplicial fans
associated to arbitrary finite root systems in recent work of S.~Fomin and
A.~Zelevinsky.
The main properties of these fans then become easy consequences of the
known
facts about
tilting modules due to K.~Bongartz, D.~Happel and C.~M.~Ringel.
\end{abstract}


\date{May 13, 2002}

\thanks{Andrei Zelevinsky's research was supported in part
by NSF grants \#DMS-9971362 and DMS-0200299. Robert Marsh's research was
supported in part by EPSRC grant GR/R17546/01.
All authors were supported in part by a University of Leicester Research
Fund Grant}

\maketitle

\tableofcontents

\section{Introduction}
\label{sec:intro}

The goal of this paper is to provide a conceptual interpretation for
\emph{generalized associahedra} introduced in \cite{FZ}.
Recall that in \cite[Section~3]{FZ}, to every finite (crystallographic)
root
system $\Phi$
is associated a smooth complete simplicial fan $\Delta (\Phi)$
in the ambient real vector space.
(Here smoothness means that every full-dimensional cone of $\Delta
(\Phi)$ is generated by a $\ZZ$-basis of the root lattice $Q$, and
completeness means that the union of all cones is the whole ambient
space.)
As proved in \cite{CFZ}, $\Delta (\Phi)$ is the normal fan of a simple
convex
polytope.
With some abuse of terminology, both $\Delta (\Phi)$ and the
corresponding simple polytope are referred to as the generalized
associahedron associated to $\Phi$.
This construction includes as special cases the well-known
associahedron (or Stasheff polytope) \cite{S} for $\Phi$ of type $A$,
and the Bott-Taubes cyclohedron \cite{BT} for $\Phi$ of type $B$ or $C$.

The construction in \cite[Section~3]{FZ} appears somewhat ad hoc (as
indicated in \cite{FZ}, it is motivated by the
developing theory of cluster algebras \cite{FZ-cl} but this
connection is still unpublished).
In this paper, we present a natural interpretation of generalized
associahedra
in terms of quiver representations.
With this interpretation, the main results in \cite[Section~3]{FZ}
become easy consequences of the known facts about tilting
representations \cite{HR}.

To explain our results in more detail, let us recall the
main ingredients of the construction in \cite{FZ}:

\begin{enumerate}

\item The set of \emph{almost positive roots}
$\Phi_{\geq -1} = \Phi_{> 0} \cup (- \Pi)$, where
$\Pi = \{\alpha_i : i \in I\}$ is a fixed set of simple roots in
the root system $\Phi$, and $\Phi_{> 0}$ is the corresponding set of
positive
roots.
Almost positive roots are precisely the generators of
one-dimensional cones in the fan $\Delta (\Phi)$.

\item The function
$\Phi_{\geq -1} \times \Phi_{\geq -1} \to \ZZ_{\geq 0}$
assigning to each $\alpha, \beta \in \Phi_{\geq -1}$
their \emph{compatibility degree} $(\alpha \| \beta)$.
A subset $C \subset \Phi_{\geq -1}$ generates a cone
in $\Delta (\Phi)$ if and only if it is \emph{compatible}, i.e.,
$(\alpha \| \beta) = 0$ for any $\alpha, \beta \in C$.

\item The group of piecewise-linear transformations of
the root lattice $Q$ generated by ``truncated simple reflections"
$\sigma_i$ for $i \in I$.
Every $\sigma_i$ preserves $\Phi_{\geq -1}$.

\item In the case when $\Phi$ is irreducible, the group generated by the
$\sigma_i$ contains two distinguished involutions $\tau_+$ and $\tau_-$
that
preserve
$\Delta (\Phi)$.

\end{enumerate}

Our aim is to provide natural interpretations of all these notions
in terms of quiver representations.
For simplicity, we only work with simply-laced root systems,
i.e., those corresponding to Dynkin diagrams of the $ADE$-type.
(The general case reduces to this with the help of a familiar
``folding"  procedure, see e.~g.~\cite[Chapter 14]{L}, or by modifying our
constructions in the
spirit of \cite{Hu,T}.)
Let $\Gamma$ be a quiver (= oriented graph) whose underlying
graph is a disjoint union
of Dynkin diagrams of type $A_n$, $D_n$, $E_6, E_7, E_8$.
By Gabriel's theorem, the dimension vector map identifies
isomorphism classes of indecomposable representations of $\Gamma$
with positive roots of the corresponding root system $\Phi$.
As a first step in our program, we include the category of
representations of $\Gamma$ into a bigger category of
\emph{decorated representations} (see Section~\ref{dqr} below).
A suitable modification of the dimension vector map establishes a
bijection $\alpha \mapsto U_\alpha^\Gamma$ between
almost positive roots in $\Phi$ and isomorphism classes of
indecomposable decorated representations of $\Gamma$.

As a next step, we find a natural interpretation for the compatibility
degree
function.
To do this, in Section~\ref{cfdqr} we construct a bifunctor $E_\Gamma$
assigning a finite-dimensional vector space to every two decorated
representations of $\Gamma$.
We then define the $\Gamma$-compatibility degree function
$\Phi_{\geq -1} \times \Phi_{\geq -1} \to \ZZ_{\geq 0}$ by setting
$(\alpha||\beta)_\Gamma=\dim E_\Gamma(U_\alpha^\Gamma,U_\beta^\Gamma)$.
As in (2) above, we call a subset $C \subset \Phi_{\geq -1}$
\emph{$\Gamma$-compatible} if $(\alpha \| \beta)_\Gamma = 0$
for any $\alpha, \beta \in C$.
We prove that the cones generated by $\Gamma$-compatible subsets
form a smooth complete simplicial fan in the
ambient space $Q_\RR$ (Theorem \ref{sf}).
This result follows without much effort from known results about
tilting representations due to D.~Happel and C.~M.~Ringel \cite{HR} and
K.~Bongartz \cite{Bo}.

The original compatibility degree and the corresponding fan $\Delta
(\Phi)$
appear as a special case of this construction corresponding to an
\emph{alternating} quiver $\Gamma_0$ for which every vertex is a source or
a
sink.
Our verification of this fact involves the final step in our
program: finding a quiver-theoretic interpretation of the
piecewise-linear involutions $\sigma_i$ acting on the root lattice
$Q$ and preserving $\Phi_{\geq -1}$.
We do this by extending the familiar reflection functors $S_i$
(see e.g., \cite{BGP}) to the functors $\Sigma_i$ on the category of
decorated
representations; more precisely, each $\Sigma_i$ sends decorated
representations of a quiver $\Gamma$ to those of the quiver $s_i \Gamma$
obtained from $\Gamma$ by reversing all arrows at a vertex $i$ which is a
source or a sink in $\Gamma$.
The main advantage of the extended reflection functors
is that each $\Sigma_i$ induces a bijection between
isomorphism classes of decorated representations of $\Gamma$ and
those of $s_i \Gamma$ (recall that the classical reflection
functors $S_i$ send some non-zero quiver representations to zero).
Furthermore, $\Sigma_i^2 (M)$ is isomorphic to $M$ for every
decorated representation $M$.
Note however that $\Sigma_i$ is \emph{not} an equivalence of
categories; in fact, it sends some non-zero morphisms to zero.

The key property of the functors $\Sigma_i$ for our purposes is
that they are compatible with the bifunctors $E_\Gamma$ in the
following sense: the vector space $E_\Gamma (M,N)$ is isomorphic to
$E_{s_i\Gamma} (\Sigma_i M, \Sigma_i N)$ for every two decorated
representations $M$ and $N$ of $\Gamma$ (Theorem \ref{ei} below).
Returning to piecewise-linear involutions $\sigma_i$, we show that
each $\sigma_i$ is obtained by combining the reflection
functor $\Sigma_i$ with the signed dimension vector map for
decorated quiver representations (see Lemma \ref{tdimsi} for a precise
statement).
This implies the following important property of the
$\Gamma$-compatibility degree:
$(\sigma_i\alpha||\sigma_i\beta)_{s_i\Gamma}=(\alpha||\beta)_\Gamma$
(Proposition \ref{sicd}).

The role of two distinguished piecewise-linear involutions
$\tau_+$ and $\tau_-$ also becomes clear in this context.
The functors $\Sigma_i$ give rise to a \emph{groupoid} whose
objects correspond to different orientations $\Gamma$ of a fixed
Dynkin graph, and each $\Sigma_i$ is viewed as a morphism from
$\Gamma$ to $s_i \Gamma$.
It is  natural to enlarge this groupoid by the \emph{duality functor} $D$
that sends every quiver $\Gamma$ to the opposite quiver 
obtained from $\Gamma$ by reversing all arrows.
We study these two groupoids in Section \ref{trg} (in a more general
setting
where the Dynkin graph is replaced by an arbitrary finite tree).
We prove in particular that for the groupoid generated by the $\Sigma_i$
and
$D$,
the automorphism group of an
alternating quiver $\Gamma_0$ is generated by two involutions
$T_+$ and $T_-$, whose induced action on the root lattice $Q$ is
given by the involutions $\tau_+$ and $\tau_-$ from \cite{FZ}.

We conclude the paper by an amusing application of our results to
the theory of ordinary quiver representations.
We say that a set $C$ of isomorphism classes of indecomposable
representations of a quiver $\Gamma$ is
\emph{Ext-free} (also called a partial tilting set) if ${\rm Ext}^1 (M,N)
= 0$ for every $M, N \in C$.
The Ext-free sets form a simplicial complex, and we show that
the $f$-vector of this complex (i.e., the number of
Ext-free sets of any given size) is an invariant of the underlying Dynkin
graph,
i.e., is independent of the choice of the orientation $\Gamma$
(Proposition \ref{pr:ext-free-enumeration}).
Note that the complexes of Ext-free sets corresponding to different
orientations of the same Dynkin graph are in general not isomorphic.
Passing from ordinary quiver representations to decorated ones
allows us to view the complex of Ext-free sets for $\Gamma$
as a subcomplex of the complex of $\Gamma$-compatible sets
discussed above; these bigger complexes are all isomorphic to each other,
by the action of extended reflection functors $\Sigma_i$.\\[3ex]
\textsc{Acknowledgments.}
This work started in July 2001 when M.~R.~and A.~Z.~visited the University
of
Leicester,
whose hospitality and financial support are gratefully acknowledged.
A.Z.~also thanks the Isaac
Newton Institute for Mathematical Sciences
in Cambridge, UK for supporting his travel to Leicester. M.~R.~ thanks
Northeastern University for kind hospitality and financial support during
his
visit in March 2002.
Finally, M.~R.~and A.~Z.~thank the organizers of the meeting
``Enveloping algebras and Algebraic Lie Representations" in April 28 - May
4,
2002,
as well as the staff at the Mathematisches Forschungsinstitut Oberwolfach,
for providing an inspiring atmosphere and the opportunity to work together
on
the finishing stages of this project.

\section{Decorated quiver representations}\label{dqr}

Let $\Gamma$ be a Dynkin quiver (that is, the underlying unoriented graph
is a
disjoint union
of Dynkin diagrams of type $A_n$, $D_n$, $E_6,E_7,E_8$), and let
$I$ be the set of vertices. Fix once and for all an algebraically closed
field
$k$. Let $\rep\Gamma$ be the category of finite
dimensional $k$-representations of $\Gamma$; the objects $M$ of
$\rep\Gamma$ thus consist of tuples of $k$-vector spaces
$(M_i)_{i\in I}$ and tuples of linear maps
$(M_a:M_i\rightarrow M_j)_{a:i\rightarrow j}$
parametrized by the arrows $a$ in $\Gamma$. For each vertex
$i\in I$, denote by $E_i$ the associated simple representation. The space
of
homomorphisms
(resp. extensions) between representations $M,N\in\rep\Gamma$
is denoted by ${\rm Hom}_\Gamma(M,N)$ (resp. ${\rm
Ext}^1_\Gamma(M,N)$). The $k$-linear dual of a $k$-vector space
$X$ is denoted by $X^*$.\\[2ex]
Let $\Phi$
be the root system associated with the Dynkin graph underlying $\Gamma$,
let $\Phi_{>0}$ be the positive roots,
and let $\Phi_{\geq-1}$ be the union of $\Phi_{>0}$ and the
negative simple roots ${-\Pi}=\{-\alpha_i\, :\, i\in I\}$. For
a subset $J\subset I$, let $\Phi(J)$ be the root subsystem of
roots whose support is contained in $J$.\\[1ex] \ Let $Q$ be the
root lattice, and let $Q_+$ be the positive span of $\Phi_{>0}$ in
$Q$. For $\gamma\in Q$ and $i\in I$, denote by $[\gamma:\alpha_i]$
the multiplicity of $\alpha_i$ in $\gamma$.\\[1ex] The
Grothendieck group of the category $\rep\Gamma$ can be
identified canonically with $Q$ via the dimension vector map
$\dimv:K_0(\rep\Gamma)\rightarrow Q$ given on objects
$M\in\rep\Gamma$ by $$\dimv(M)=\sum_{i\in
I}(\dim M_i)\alpha_i.$$

Associate to $\Gamma$ a new quiver
$\widetilde{\Gamma}=\Gamma\cup\{i^-\, :\, i\in I\}$, with no
arrows pointing from or to the added vertices $i^-$. The category
$\drc$ will be called the decorated representation category for $\Gamma$.
The
simple representation $E_{i^-}$ will be denoted by $E_i^-$; the
objects of $\drc$ are thus of the form $M=M^+\oplus V$, where
$M^+\in\rep\Gamma$, and $V=\bigoplus_{i\in I}(E_i^-\otimes V_i)$
is viewed as an $I$-graded vector space.
More explicitly, an object of $\drc$ consists of tuples of vector spaces
$M_i$
and $V_i$
for $i\in I$, together with linear maps $M_a:M_i\rightarrow M_j$ for all
arrows
$a:i\rightarrow j$ in $\Gamma$; morphisms between two such pairs of tuples
$(M_i,V_i)$
and $(N_i,W_i)$ consist of linear maps $f_i:M_i\rightarrow N_i$,
together with linear maps $\varphi_i:V_i\rightarrow W_i$,
such that $f_j\circ M_a=N_a\circ f_i$ for all arrows $a:i\rightarrow
j$.\\[1ex]
The space of $I$-graded
homomorphisms of $I$-graded vector spaces $V,W$ will be
denoted by ${\rm Hom}^I(V,W)$.\\[1ex]
The category $\rep\Gamma$
is viewed as a subcategory of $\drc$ in the obvious way. Define
the signed dimension vector $\tdimv(M)$ of $M=M^+\oplus V$ by
$$\tdimv(M)=\dimv M^+-\sum_{i\in I}(\dim V_i)\alpha_i=\sum_{i\in I}(\dim
M_i-\dim V_i)\alpha_i.$$

By Gabriel's theorem, the dimension vector induces a bijection
between the isoclasses of indecomposable representations in
$\rep\Gamma$ and the positive roots $\Phi_{>0}$.\\[2ex] We have
the following extension of Gabriel's theorem, which is immediate
from the definition of the decorated representation category:

\begin{proposition} The map $\tdimv$ induces a bijection
$\alpha\mapsto[U_\alpha^\Gamma]$
between the roots $\alpha\in\Phi_{\geq -1}$ and isoclasses
$[U_\alpha^\Gamma]$ of indecomposable objects in $\drc$.
\end{proposition}

\section{Clusters from decorated representations}\label{cfdqr}

For decorated representations $M=M^+\oplus V$ and $N=N^+\oplus W$ in
$\drc$,
we define
\begin{equation}\label{def:egamma}
E_\Gamma(M,N)={\rm Ext}^1_\Gamma(M^+,N^+)\oplus{\rm
Ext}^1_\Gamma(N^+,M^+)^*\oplus\end{equation}
$${\rm Hom}^I(M^+,W)\oplus{\rm
Hom}^I(V,N^+).
$$
Obviously, this defines a bifunctor
$$E_\Gamma(\_,\_):(\drc)^{\rm op}\times\drc\rightarrow\vect,$$
contravariant in the first argument, and covariant in the second one. Here
$\vect$ denotes the category of finite-dimensional vector spaces.

In view of (\ref{def:egamma}), we have the following
characterization of decorated representations $M=M^+\oplus V$
such that $E_\Gamma(M,M)=0$:
\begin{eqnarray}
\label{eq:E-Gamma=0}
&&\text{$E_\Gamma(M,M)=0$ if and only if ${\rm
Ext}^1_\Gamma(M^+,M^+)=0$,}\\
\nonumber &&\text{and $M^+$ and $V$ have disjoint
supports in $I$.}
\end{eqnarray}

We define the $\Gamma$-compatibility degree of two roots
$\alpha,\beta\in\Phi_{\geq -1}$ by
\begin{equation}\label{compdegree}
(\alpha||\beta)_\Gamma=\dim
E_\Gamma(U_\alpha^\Gamma,U_\beta^\Gamma).
\end{equation}

This definition makes it clear that $(||)_{\Gamma}$ is symmetric and
compatible
with root subsystems in the
following sense: if $\Gamma(J)$ denotes the full subquiver of
$\Gamma$ supported on a subset of vertices $J\subset I$, and if
$\alpha,\beta\in\Phi_{\geq-1}(J)$, then
$(\alpha||\beta)_{\Gamma(J)}=(\alpha||\beta)_\Gamma$ (compare
\cite[Proposition 3.~3.~]{FZ}).
There is an additional symmetry property involving the natural duality
functor
$D:\rep\widetilde{\Gamma}\rightarrow\rep\widetilde{\Gamma}^{\rm op}$ that
acts
by
replacing all vector spaces and maps by their duals, thus replacing
$\Gamma$
with the opposite orientation $\Gamma^{\rm op}$:

\begin{proposition}\label{eiandd} For all roots
$\alpha,\beta\in\Phi_{\geq-1}$,
we have $(\alpha||\beta)_\Gamma=(\alpha||\beta)_{\Gamma^{\rm op}}$.
\end{proposition}

This follows at once from the natural identification
$$E_{\Gamma^{\rm op}}(DN,DM)\simeq E_\Gamma(M,N).$$

We can now define
$\Gamma$-compatible subsets and $\Gamma$-clusters as in \cite[Definition
3.~4.~]{FZ}.

\begin{definition}\label{def:compatible} {\rm A subset $C$ of $\Phi_{\geq
-1}$
is called \emph{
$\Gamma$-compatible} if
$(\alpha||\beta)_{\Gamma}=0$ for all
$\alpha,\beta\in C$. The subset $C$ is called a \emph{$\Gamma$-cluster} if
it
is a
maximal (by inclusion) $\Gamma$-compatible subset.}
\end{definition}

\begin{definition} {\rm The \emph{negative support}
$S(C)$ of a subset $C$ of $\Phi_{\leq-1}$ is defined by $S(C)=\{i\in I\,
:\, -\alpha_i\in C\}$. The subset $C$ is called
\emph{positive} if $C\subset\Phi_{>0}$, i.e., $S(C)=\emptyset$.}
\end{definition}

Moreover, \cite[Proposition 3.~6.~]{FZ} holds with the same
proof.

\begin{proposition}\label{redpos} For every subset $J\subset I$, the
correspondence $C\mapsto
C-\{-\alpha_i\, :\, i\in J\}$ defines a bijection between
$\Gamma$-compatible subsets (resp.~clusters) for $\Phi$ with negative
support $J$ and positive $\Gamma$-compatible subsets (resp.~clusters) for
$\Phi(I- J)$.
\end{proposition}

We can now prove the main properties of $\Gamma$-clusters.

\begin{proposition}[Purity]\label{purity} All $\Gamma$-clusters have the
same
size $n$, the rank of
the root system
$\Phi$. Moreover, each $\Gamma$-cluster is a $\ZZ$-basis for the root
lattice $Q$.
\end{proposition}

\proof In view of Proposition \ref{redpos}, it suffices to
prove our statements for positive $\Gamma$-clusters.
To any subset $C \subset \Phi_{> 0}$ let us associate a
$\Gamma$-representation
$X_C=\bigoplus_{\alpha\in C} U_\alpha^\Gamma$.
By our definitions, $C$ is $\Gamma$-compatible if and only if ${\rm
Ext}^1(X_C,X_C)=0$.
Now we use the theory of tilting modules introduced in \cite{HR}.
In our situation, they can be defined as follows (see \cite[Theorem
4.5]{HR}):
a $\Gamma$-representation $T$ is \emph{tilting} if ${\rm Ext}^1(T,T)=0$,
and
$T$
has precisely $n$ different isoclasses of indecomposable direct summands.

Now let $C$ be a positive $\Gamma$-cluster.
By \cite[Lemma 2.1]{Bo}, $X_C$ is a direct summand of
a tilting module $T$.
By the maximality property of $C$, the representation $X_C$ contains all
isoclasses of indecomposable direct summands of $T$
and so $|C| = n$. The fact that $C$ is a $\ZZ$-basis for $Q$ is contained
in
the proof of \cite[Lemma 4.3]{HR}. \endproof

\begin{proposition}[Cluster expansion]
\label{clexp} Each element $\gamma\in Q$
has a unique
$\Gamma$-cluster expansion,
that is, $\gamma=\sum_{\alpha\in\Phi_{\geq-1}}m_\alpha\alpha$ such
that $m_\alpha\in\ZZ_{\geq 0}$ for all $\alpha$, and $m_\alpha m_\beta=0$
whenever $\alpha,\beta\in\Phi_{\geq-1}$
are not $\Gamma$-compatible.
\end{proposition}

\proof We prove that for each $\gamma\in Q$, there exists a unique
decorated representation $M=M^+\oplus V\in\drc$ such that $\tdimv
M=\gamma$ and
$E_\Gamma(M,M)=0$ (the proposition follows from this by definition
of the $\Gamma$-compatibility degree).
In view of (\ref{eq:E-Gamma=0}),
the condition $E_\Gamma(M,M)=0$ means that ${\rm
Ext}^1_\Gamma(M^+,M^+)=0$, and $M^+$ and $V$ have disjoint
supports in $I$.
Denote by $d=\sum_{i\in I}d_i\alpha_i$ the dimension vector of $M^+$.
Note that it is uniquely determined by $\gamma$, namely
$d=\sum_{i\in I}\max([\gamma:~\alpha_i],0)\alpha_i$.
Define the variety $R_d$ of $k$-representations of $\Gamma$ of
dimension vector $d$ by $R_d=\oplus_{a:i\rightarrow j}{\rm
Hom}(k^{d_i},k^{d_j})$.
The group $G_d=\prod_{i\in I}{\rm GL}(k^{d_i})$ acts on $R_d$ via
$(g_i)(M_a)=(g_jM_a g_i^{-1})_{a:i\rightarrow j}$.
By definition, the orbits of $G_d$ in $R_d$ correspond bijectively to the
isomorphism classes of
representations of $\Gamma$ of dimension vector $d$.
Denote by ${\mathcal O}_A$ the $G_d$-orbit of a representation $A$.
Since $R_d$ is irreducible and there are only finitely many orbits by
Gabriel's theorem,
$R_d$ contains a unique dense orbit. It is well known that (see, e.g.,
\cite{Bo-geo})
$$\codim_{R_d}{\mathcal O}_A=\dim{\rm
Ext}^1_{\Gamma}(A,A).$$
It follows that there exists a unique (up to
isomorphism) representation $M^+$ of dimension vector $d$
such that ${\rm Ext}^1_\Gamma(M^+,M^+)=0$. Since the $I$-graded
$k$-space $V$ is uniquely determined by its signed dimension vector
$\gamma-d$, uniqueness and existence follow.\endproof

As in \cite[Proof of
Theorem 1.~10.~]{FZ}, Propositions \ref{purity} and \ref{clexp} imply the
following:

\begin{theorem}\label{sf} The simplicial cones ${\bf R}_{\geq 0}C$
generated
by all
$\Gamma$-clusters $C$ form a smooth complete
simplicial fan $\Delta_\Gamma$ in the ambient real vector space $Q_{\bf
R}$ of
the
root system.
\end{theorem}

\begin{remark} \label{ungerremark}
{\rm The ``positive part" of $\Delta_\Gamma$ was studied by C.~Riedt\-mann
and A.~Scho\-field \cite{RS} and L.~Unger \cite{U}
following a suggestion of C.~M.~Ringel.}
\end{remark}

\section{Extended reflection functors}\label{reflectionfunctors}

We recall the well-known reflection functors on quiver
representations \cite{BGP}. Let $i\in I$ be a source (resp.~sink) in
$\Gamma$, and let $s_i\Gamma$ be the quiver obtained by reversing
all arrows pointing from (resp.~to) $i$. Given a representation
$M\in\rep\Gamma$, define a representation
$S_i(M)\in\rep s_i\Gamma$ by
\begin{equation}\label{def:classicalsi}
S_i(M)_i=\left\{\begin{array}{lcl} {\rm
Coker}(\bigoplus_{a:i\rightarrow j}M_a)&,&\mbox{$i$ a
source}\\ {\rm Ker}(\bigoplus_{a:j\rightarrow
i}M_a)&,&\mbox{$i$ a sink}
\end{array}\right\},\;\; S_i(M)_j=M_j\mbox{ for all }j\not=i,
\end{equation}
and the maps $S_i(M)_a$ are either unchanged, or derived in an obvious way
from
the canonical injection or projection maps. The functor $S_i$ is defined
on
morphisms using the canonical induced maps on cokernels (resp.
kernels).\\[1ex]
Let $\rep^i\Gamma$ denote the full
subcategory of $\rep\Gamma$ consisting of representations without $E_i$ as
a direct summand. Then $S_i$ induces an equivalence of categories
$$S_i:\rep^i\Gamma\stackrel{\sim}{\rightarrow}\rep^is_i\Gamma.$$

\subsection{Definition and basic properties}\label{rf}

Assume that $i\in I$ is a source or a sink in $\Gamma$. We will
now define a reflection functor
$$\Sigma_i:\drc\rightarrow\rep\widetilde{s_i\Gamma},$$ which has
the advantage -- over the usual reflection functor $S_i$ -- that it
induces
a bijection on isomorphism classes of decorated representations.\\[1ex]
Let
$M=M^+\oplus V$ be an object of $\drc$. Suppose first that $i$ is
a source in $\Gamma$, and consider the canonical exact sequence of
$\Gamma$-representations
\begin{equation}\label{seqm+m'}
\ses{E_i\otimes{\rm Hom}_\Gamma(E_i,M^+)}{M^+}{M'}.
\end{equation}
Note for future reference the following well-known exact sequences:
\begin{equation}\label{useful}
0\rightarrow{\rm Hom}_\Gamma(E_i,M^+)\rightarrow
M_i\rightarrow \bigoplus_{i\rightarrow j}M_j\rightarrow {\rm
Ext}^1_\Gamma(E_i,M^+)\rightarrow 0.
\end{equation}
The analogue for $i$ a sink is
\begin{equation}\label{useful*}
0\rightarrow {\rm Hom}_\Gamma(M^+,E_i)\rightarrow
M_i^*\rightarrow\bigoplus_{j\rightarrow i}M_j^*
\rightarrow{\rm Ext}^1_\Gamma(M^+,E_i)\rightarrow 0.
\end{equation}
Both these sequences are special cases of \cite[Lemma 2.1]{Ri}.

In particular, (\ref{useful}) implies the natural isomorphism
$${\rm Hom}_\Gamma(E_i,M^+)\simeq {\rm Ker}(\bigoplus_{a:i\rightarrow
j}M_a).$$
Returning to (\ref{seqm+m'}), it is clear that ${\rm
Hom}_\Gamma(E_i,M')=0$,
thus
$M'\in\rep^i\Gamma$. We define $\Sigma_i(M)=\overline{M}^+\oplus
\overline{V}$
by setting
\begin{equation}\label{defsi}
\overline{M}^+=S_i(M')\oplus
(E_i\otimes V_i),\quad
\overline{V}=(E_i^-\otimes {\rm
Hom}_\Gamma(E_i,M^+))\oplus\bigoplus_{j\not= i}(E_j^-\otimes
V_j).
\end{equation}
More explicitly, $\Sigma_i(M)$ replaces $M_i$ by ${\rm
Coker}(\bigoplus_{a:i\rightarrow j}M_a)\oplus V_i$, replaces $V_i$ by
${\rm
Ker}(\bigoplus_{a:i\rightarrow j}M_a)$,
and leaves the rest of the spaces unchanged.\\[2ex]
To define $\Sigma_i$ on morphisms, let $N=N^+\oplus W$ be
another object of $\drc$, and define $N'$, $\overline{N}$, $\overline{W}$
in
the same way as above.
Let $f=f^+\oplus(\bigoplus_j\varphi_j)$
be a morphism from $M$ to $N$ in $\drc$. Since ${\rm
Hom}_\Gamma(E_i,N')=0$, the map $f^+$ induces maps
$$\psi:{\rm Hom}_\Gamma(E_i,M^+)\rightarrow{\rm Hom}_\Gamma(E_i,N^+),\quad
f':M'\rightarrow N'.$$
We then define the morphism $\Sigma_i(f)$ as the direct sum of the
morphisms
$$S_i(f')\oplus\varphi_i:\overline{M}^+\rightarrow\overline{N}^+$$
and
$$\psi\oplus\bigoplus_{j\not=i}\varphi_j:\overline{V}\rightarrow\overline{W}.$$
More explicitly, $\Sigma_i(f)$ replaces the map $f_i$ by $g\oplus
\varphi_i$
(with respect to the above sum decomposition),
where $g$ is the map induced by $f_i$ on the cokernels,
replaces $\varphi_i$ by the map $h$ induced on the kernels, and leaves the
rest of the maps unchanged.

To treat the case when $i$ is a sink, we use the duality functor $D$ (see
Section \ref{cfdqr}).
This allows us to reduce the case of a sink to that of a source by
defining
$$\Sigma_i=D\circ\Sigma_i\circ D.$$
Explicitly, take the canonical exact sequence
$$\ses{M'}{M^+}{E_i\otimes {\rm Hom}_\Gamma(M^+,E_i)^*}.$$
Then we have $\Sigma_i(M)=\overline{M}^+\oplus\overline{V}$, where
$$\overline{M}^+=S_i(M')\oplus(E_i\otimes V_i)$$
and
$$\overline{V}=(E_i^-\otimes {\rm
Hom}_\Gamma(M^+,E_i)^*)\oplus\bigoplus_{j\not=i}(E_j^-\otimes
V_j).$$
The description of the action of $\Sigma_i$ on morphisms is similar to the
case of a source and is left to the reader.

The following proposition summarizes some simple properties of the
extended
reflection functors.
All of them are immediate from the definitions.
\begin{proposition}\label{pr:easycons}
\begin{enumerate}
\item If $M\in\rep^i\Gamma\subset\drc$, then $\Sigma_i(M)=S_i(M)$.
\item We have $\Sigma_i(E_i)=E_i^-$ and
$\Sigma_i(E_i^-)=E_i$.
\item We have $\Sigma_i(E_j^-)=E_j^-$ whenever $i\not=j$.
\item Suppose $i$ is a source (resp.~a sink) in $\Gamma$, and
$M\in\rep^i\Gamma\subset\drc$.
Let $f$ be any morphism from $M$ to $E_i$ (resp.~from $E_i$ to $M$). Then
$\Sigma_i(f)=0$.
\end{enumerate}
\end{proposition}

\begin{proposition}\label{pr:si2} For any decorated representation $M\in
\drc$
and any vertex
$i\in I$ we have $\Sigma_i^2(M)\simeq M$.
\end{proposition}

\proof Assume $i\in I$ to be a source in $\Gamma$ (the case of a sink
follows
by duality).
For a decorated representation $M=M^+\oplus V\in\drc$, denote by
$N=N^+\oplus
W$
the decorated representation $\Sigma_i(M)$. (In the notation of
\ref{defsi},
we have $N=\overline{M}$ and $W=\overline{V}$.)
By definition, we have:
$$\Sigma_i^2(M)=\overline{N}^+\oplus\overline{W}.$$
Here
$$\overline{N}^+={S_i(N')\oplus (E_i\otimes W_i)}\simeq{M'\oplus
(E_i\otimes{\rm Hom}_\Gamma(E_i,M^+)},$$
and
\begin{eqnarray*}
\overline{W}&=&{(E_i^-\otimes{\rm Hom}_{s_i\Gamma}(S_i(N')\oplus
(E_i\otimes
V_i),E_i))^*\oplus\bigoplus_{j\not=i}(E_j^-\otimes W_j)}\\
&\simeq&{(E_i^-\otimes V_i)\oplus\bigoplus_{j\not=i}(E_j^-\otimes
V_j)}\simeq
V.
\end{eqnarray*}
Note that ${M'\oplus (E_i\otimes{\rm Hom}_\Gamma(E_i,M^+)}$ is
(non-canonically) isomorphic to $M^+$
since the sequence (\ref{seqm+m'}) splits. Therefore $\Sigma_i^2(M)\simeq
M^+\oplus V=M$. \endproof

\begin{corollary}\label{cor:isogroth} The functor $\Sigma_i$ induces an
isomorphism of
Gro\-then\-dieck groups
$\Sigma_i:K_0(\drc)\stackrel{\sim}{\rightarrow}K_0(\rep
\widetilde{s_i\Gamma})$.
\end{corollary}

We define piecewise linear reflections $\sigma_i$ on $Q$ as in
\cite[2.2]{FZ}.
For the simply-laced case,
this can be stated as follows:

\begin{definition}\label{def:sigmai} For $i\in I$, define
$\sigma_i:Q\rightarrow Q$ by
$$[\sigma_i\alpha:\alpha_j]=\left\{\begin{array}{ccc}
[\alpha:\alpha_j]&,&j\not=i\\[1ex]
-[\alpha:\alpha_i]+\displaystyle\sum_{k-\!\!\!-\!\!\!-
i}\max([\alpha:\alpha_k],0)&,&j=i\end{array}\right.,$$
where $k-\!\!\!-\!\!\!- i$ means that $k$ and $i$ are joined by an edge in
$I$.
\end{definition}

\begin{lemma}\label{tdimsi} For any source or sink $i\in I$ of $\Gamma$,
and
any object
$M=M^+\oplus V\in\drc$ such that $\dim M_i\cdot\dim
V_i=0$, we have $\tdimv(\Sigma_i(M))=\sigma_i(\tdimv(M))$.
\end{lemma}

\proof Let $i$ be a source in $\Gamma$ (the case of a sink is dual).
The reflection functor $\Sigma_i$ does not change any component
$M_j$ or $V_j$ for $j\not=i$, which corresponds to the first
part of the above definition of $\sigma_i$. So we concentrate on
the $i$-components. By definition,
we have $[\tdimv M:\alpha_i]=\dim M_i-\dim V_i$. Using the
definition of $\Sigma_i$, we compute:
$$[\tdimv\Sigma_i(M):\alpha_i]=\dim S_i(M')_i+\dim V_i-\dim{\rm
Hom}_\Gamma(E_i,M).$$
In view of the exact sequence (\ref{useful}),
$$\dim S_i(M')_i-\dim{\rm Hom}_\Gamma(E_i,M)=-\dim M_i+\sum_{i\rightarrow
j}\dim M_j.$$
Therefore,
\begin{eqnarray*}
[\tdimv\Sigma_i(M):\alpha_i]&=&-\dim M_i+\sum_{i\rightarrow j}\dim
M_j+\dim V_i=\\ &=&-[\tdimv M:\alpha_i]+\sum_{i\rightarrow
j}\dim M_j.
\end{eqnarray*}
Under the assumption $\dim M_i\cdot\dim V_i=0$, this last
expression is the same as in Definition \ref{def:sigmai}. \endproof

\begin{corollary} \label{correspondence}
For all roots $\alpha\in\Phi_{\geq -1}$, and any $i\in I$ which is a
source or
a sink in $\Gamma$, we
have
$$\tdimv(\Sigma_i(U_\alpha^\Gamma))=\sigma_i(\tdimv(U_\alpha^\Gamma)).$$
\end{corollary}

\proof Notice that for any indecomposable object
$M\in\drc$, we have $\dim M_i\cdot\dim V_i=0$ for all $i\in I$.
\endproof

\subsection{Extended reflection functors and compatibility
degree}\label{racd}


Recall the bifunctor $E_\Gamma$ on decorated representations introduced in
(\ref{def:egamma}).

\begin{theorem}\label{ei} If $i\in I$ is a source or a sink in $\Gamma$
then
$E_{s_i\Gamma}(\Sigma_i(M),\Sigma_i(N))\simeq
E_\Gamma(M,N)$ for any decorated representations $M$ and $N$.
\end{theorem}

\proof We start with a lemma which is an immediate consequence of
(\ref{useful}) and (\ref{useful*}).

\begin{lemma}\label{l:homintsi} Let $A$ be a representation of $\Gamma$.
If
$i$ is a source in $\Gamma$, then
$$S_i(A)_i\simeq{\rm Ext}_\Gamma^1(E_i,A).$$
If $i$ is a sink in $\Gamma$, then
$$S_i(A)_i\simeq{\rm Ext}^1_\Gamma(A,E_i)^*.$$
\end{lemma}

Returning to the proof of Theorem \ref{ei}, we will only treat the case of
a
source $i\in I$; the case of a
sink is dual. As usual, we write $M=M^+\oplus V$ and $N=N^+\oplus
W$. Recall the (non-canonical) decompositions
$$M^+=M'\oplus (E_i\otimes X),\quad N^+=N'\oplus (E_i\otimes Y),$$
where $X={\rm Hom}_\Gamma(E_i,M^+)$ and $Y={\rm Hom}_\Gamma(E_i,N^+)$.
Using the
definition of $\Sigma_i$, we decompose\\
$E_{s_i\Gamma}(\Sigma_i(M),\Sigma_i(N))$
\begin{eqnarray*}
&=&{\rm Ext}^1_{s_i\Gamma}(S_i(M'),S_i(N'))
\oplus({\rm Ext}^1_{s_i\Gamma}(S_i(M'),E_i)\otimes W_i)
\oplus\\
&&{\rm Ext}^1_{s_i\Gamma}(S_i(N'),S_i(M'))^*
\oplus ({\rm Ext}^1_{s_i\Gamma}(S_i(N'),E_i)\otimes V_i)^*
\oplus\\
&&{\rm Hom}(S_i(M')_i,Y)\oplus{\rm
Hom}(V_i,Y)\oplus\bigoplus_{j\not=i}{\rm Hom}(M_j,W_j)\oplus\\
&&{\rm Hom}(X,S_i(N')_i)\oplus{\rm
Hom}(X,W_i)\oplus\bigoplus_{j\not=i}{\rm Hom}(V_j,N_j).
\end{eqnarray*}
Since the classical reflection functor $S_i$ induces an equivalence of
categories
$\rep^i\Gamma\stackrel{\sim}{\rightarrow}\rep^is_i\Gamma$, we have
\begin{eqnarray}\label{extsi}
&&{\rm
Ext}^1_{s_i\Gamma}(S_i(M'),S_i(N'))\simeq {\rm Ext}^1_\Gamma(M',N'),\\
&&{\rm Ext}^1_{s_i\Gamma}(S_i(N'),S_i(M'))^*{\simeq {\rm
Ext}^1_\Gamma(N',M')^*}.
\end{eqnarray}
Furthermore, applying the second part of Lemma \ref{l:homintsi} to
$A=S_i(M')$
and using
Proposition \ref{pr:si2}, we obtain
\begin{equation}\label{homintsi}
{\rm Ext}^1_{s_i\Gamma}(S_i(M'),E_i)\otimes W_i{\simeq {\rm
Hom}(M_i',W_i)}.
\end{equation}
Similarly
\begin{equation}\label{homintsi*}
({\rm
Ext}^1_{s_i\Gamma}(S_i(N'),E_i)\otimes V_i)^*{\simeq{\rm Hom}(V_i,N_i')}.
\end{equation}
Substituting the expressions provided by (\ref{extsi}) - (\ref{homintsi*})
into the above
decomposition of $E_{s_i\Gamma}(\Sigma_i(M),\Sigma_i(N))$ and regrouping
the
terms,
we can write it as the direct sum of the following four terms:
\begin{eqnarray}
\label{term1}&&({\rm Ext}_\Gamma^1(M',N')\oplus{\rm
Hom}(X,S_i(N')_i))\oplus\\
\label{term2}&&({\rm Ext}_\Gamma^1(N',M')^*\oplus{\rm
Hom}(S_i(M')_i,Y))\oplus\\
\label{term3}&&({\rm Hom}(M_i',W_i)\oplus\bigoplus_{j\not=i}{\rm
Hom}(M_j,W_j)\oplus{\rm Hom}(X,W_i))\oplus\\
\label{term4}&&({\rm Hom}(V_i,N_i')\oplus{\rm
Hom}(V_i,Y)\oplus\bigoplus_{j\not=i}{\rm Hom}(V_j,N_j)).
\end{eqnarray}

To complete the proof, we show that these four terms can be identified
with
the four terms in (\ref{def:egamma}) in the same order.

First we have
\begin{eqnarray*}
{\rm Ext}^1_\Gamma(M^+,N^+)&=&{\rm Ext}^1_\Gamma(M'\oplus
(E_i\otimes X), N'\oplus (E_i\otimes Y))\\ &=&{\rm
Ext}^1_\Gamma(M',N')\oplus({\rm Ext}^1_\Gamma(E_i,N')\otimes X^*)
\end{eqnarray*}

This agrees with (\ref{term1}) by Lemma \ref{l:homintsi}.
The identification of ${\rm Ext}^1_\Gamma(N^+,M^+)^*$ with (\ref{term2})
is
proved dually.
Finally, the terms in (\ref{term3}) and (\ref{term4}) obviously equal the
last
two terms in (\ref{def:egamma}). \endproof

\begin{remark} {\rm The isomorphism given in the theorem is not natural,
since
the $\Sigma_i$ do not provide an equivalence of categories. But it
can be shown (with the same proof as above) that there exists an
isomorphism of functors $E_{s_i\Gamma}(\Sigma_i,\Sigma_i)\simeq
E_\Gamma(\Sigma_i^2,\Sigma_i^2)$.}
\end{remark}

As an immediate consequence of the definition (\ref{compdegree}),
Corollary \ref{correspondence} and Theorem \ref{ei}, we get:

\begin{proposition}\label{sicd} If $i\in I$ is a source or a sink in
$\Gamma$,
then $(\sigma_i\alpha||\sigma_i\beta)_{s_i\Gamma}=(\alpha||\beta)_\Gamma$
for
all $\alpha,\beta\in\Phi_{\geq -1}$.
\end{proposition}

This result has the following implications.

\begin{corollary}\label{cor:independence} If $\Gamma$ and $\Gamma'$ are
two
quivers with the same underlying Dynkin graph,
then the simplicial complexes $\Delta_\Gamma$ and $\Delta_{\Gamma'}$ are
isomorphic to each other.
\end{corollary}

\proof If $\Gamma'=s_i\Gamma$ for some source or sink $i$, then by
Proposition
\ref{sicd},
the permutation $\sigma_i$ of $\Phi_{\geq -1}$ transforms $\Delta_\Gamma$
to
$\Delta_{\Gamma'}$.
Any two quivers can be obtained from each other by a sequence of such
transformations
(see, e.~g.~ \cite[Theorem 1.2]{BGP}). \endproof

We now show that the notion of $\Gamma$-compatibility degree
includes as a special case the compatibility degree introduced in
\cite[Section~3.1]{FZ}.
Let $\Gamma_0$ be an alternating quiver (i.~e.~such that each
vertex is a source or a sink).
Denote the set of sources of $\Gamma_0$ by $I^+$, and the set of sinks by
$I^-$.
Following \cite{FZ}, we define the transformations $\tau_+$ and
$\tau_-$ by
$$\tau_+ = \prod_{i\in I^+}\sigma_i, \,\,
\tau_- = \prod_{i\in I^-}\sigma_i.$$
Recall from \cite{FZ} that the compatibility
degree is characterized by the following two
properties:
\begin{eqnarray}
\label{eq:comp-1}
&&\text{$(-\alpha_i||\beta)=\max([\beta:\alpha_i],0)$ for all $i\in I$ and
$\beta\in\Phi_{\geq-1}$;}\\
\label{eq:comp-2}
&&\text{$(\tau_\epsilon\alpha||\tau_\epsilon\beta)=(\alpha||\beta)$
for $\epsilon\in\{+,-\}$.}
\end{eqnarray}

\begin{corollary}\label{cfz} The $\Gamma_0$-compatibility degree coincides
with the compatibility degree defined by {\rm
(\ref{eq:comp-1})-(\ref{eq:comp-2})}.
Therefore, $\Gamma_0$-clusters coincide with the clusters in the sense
of~\cite{FZ}.
\end{corollary}

\proof
It suffices to show that the $\Gamma_0$-compatibility degree satisfies
(\ref{eq:comp-1})-(\ref{eq:comp-2}).
The property (\ref{eq:comp-1}) follows at once from
(\ref{compdegree}).
For (\ref{eq:comp-2}), notice that $\tau_+$ and $\tau_-$ send $\Gamma_0$
to
$\Gamma_0^{\rm op}$, and apply Propositions \ref{sicd} and \ref{eiandd}.
\endproof

\section{The reflection groupoid}\label{trg}

In this section, we show that the
group generated by the two involutions $\tau_+,\tau_-$ has
a natural interpretation in the context of decorated quiver
representations.
To do this, we consider the following more general setup.
Suppose that $I$ is an arbitrary finite tree. We associate to $I$ two
groupoids defined as follows.

Let $R_0(I)$ denote the groupoid whose objects are given by the set ${\rm
Quiv}(I)$
of all quivers obtained by specifying an orientation of $I$, and whose
morphisms are defined as follows.
First of all, for any $\Gamma\in {\rm Quiv}(I)$ and any source or sink $i$
in
$\Gamma$, there is an ``elementary" isomorphism
$\Sigma_i:\Gamma\stackrel{\sim}{\rightarrow}s_i\Gamma$.
The morphisms in $R_0(I)$ are arbitrary compositions of these
isomorphisms, subject to the following relations
(where $1$ denotes the appropriate identity automorphism):
\begin{itemize}
\item[(R1)] $\Sigma_i^2=1$,
\item[(R2)] $\Sigma_i\Sigma_j=\Sigma_j\Sigma_i$ whenever $i$ and $j$ are
not linked in $I$.
\end{itemize}
In both cases, the relation holds whenever the maps and their compositions
are defined.

We also consider another groupoid $R(I)$ with the same set of objects
${\rm Quiv}(I)$, obtained from $R_0 (I)$ by adjoining one more
``elementary
isomorphism"
\mbox{$D:\Gamma\stackrel{\sim}{\rightarrow}\Gamma^{\rm op}$},
for any $\Gamma\in {\rm Quiv}(I)$, satisfying the following relations
(as before, they must hold whenever the morphisms and compositions on both
sides are defined):
\begin{itemize}
\item[(R3)] $D^2=1$,
\item[(R4)] $D\Sigma_i=\Sigma_iD$ for all $i\in I$.
\end{itemize}
Let $I=I^+\cup I^-$ be a disjoint union of $I$ into two completely
disconnected subsets, i.e. no two vertices in $I^+$ (respectively, $I^-$)
have a common edge. Such a decomposition exists because $I$ is a tree,
and is unique up to interchanging $I^+$ and $I^-$.

Denote by $\Gamma_0$ the alternating orientation of $I$ in which every
element of $I^+$ is a source.
Let $\Sigma_+=\prod_{i\in I^+}\Sigma_i$ and $\Sigma_-=\prod_{i\in
I^-}\Sigma_i$.
Note that, by the relation (R2), these morphisms
do not depend on the order in which the compositions are taken.
We can regard $\Sigma_+,\Sigma_-$ as morphisms (in either of the groupoids
$R_0(I)$ or $R(I)$) from $\Gamma_0$ to the opposite orientation
$\Gamma_0^{\rm op}$, or from $\Gamma_0^{\rm op}$ to $\Gamma_0$.
The main result of this section can now be stated as follows:

\begin{theorem} \label{alternatingaut}
(a) The automorphism group ${\rm Aut}_{R_0(I)}(\Gamma_0)$ is generated by
$\Sigma_+ \Sigma_-$ and $\Sigma_- \Sigma_+$. \\
(b) The automorphism group ${\rm Aut}_{R(I)}(\Gamma_0)$ is generated by
$D\Sigma_+$ and $D\Sigma_-$.
\end{theorem}

Before embarking on the proof of this theorem, we note the following:

\begin{remark} \rm
Let $T_+=D\Sigma_+:\Gamma_0\stackrel{\sim}{\rightarrow}\Gamma_0$ and
let $T_-=D\Sigma_-:\Gamma_0\stackrel{\sim}{\rightarrow}\Gamma_0$.
Theorem~\ref{alternatingaut} states that $T_+$ and $T_-$ generate the
automorphism group ${\rm Aut}_{R(I)}(\Gamma_0)$. In the case where
$I$ is a simply-laced Dynkin diagram, we can interpret the
isomorphisms $\Sigma_i$ in $R(I)$ as the extended reflection functors of
Section~\ref{reflectionfunctors} and the isomorphism $D$ as the usual
duality ${\rep\widetilde{\Gamma}}\stackrel{\sim}{\rightarrow}
{\rep\widetilde{\Gamma^{\rm op}}}$.
Using Corollary~\ref{correspondence} and the fact that
$\tdimv(DM)=\tdimv(M)$ for all objects $M$ of $\drc$ (for any
$\Gamma$), we see that if
$T_+$ and $T_-$ are regarded as bijections of the set of objects of
${\rep\widetilde{\Gamma}_0}$, they induce the involutions $\tau_+$ and
$\tau_-$ of the root lattice $Q$ introduced in~\cite{FZ}.
\end{remark}

Theorem~\ref{alternatingaut} has the following corollary:
\begin{corollary}
Let $\Gamma$ be any element of ${\rm Quiv}(I)$. Then: \\
(a) The automorphism group ${\rm Aut}_{R_0(I)}(\Gamma)$ is cyclic. \\
(b) The automorphism group ${\rm Aut}_{R(I)}(\Gamma)$ is generated by
two involutions.
\end{corollary}

\proof This follows immediately from Theorem~\ref{alternatingaut}, using
the fact that ${\rm Quiv}(I)$ is connected under the operations $s_i$
(by \cite[Theorem 1.2]{BGP}), and so all automorphism groups in the
groupoid
are conjugate.
\endproof

\noindent {\bf Proof of Theorem~\ref{alternatingaut}.}
We note first of all that (b) is an easy formal consequence of (a):
an element of ${\rm Aut}_{R_0(I)}(\Gamma_0)$ is either
$(\Sigma_+ \Sigma_-)^k = (D \Sigma_+  D \Sigma_-)^k$, or
$(\Sigma_- \Sigma_+)^k = (D \Sigma_-  D \Sigma_+)^k$;
and an element of ${\rm Aut}_{R(I)} (\Gamma_0) \setminus {\rm
Aut}_{R_0(I)}
(\Gamma_0)$ can be written as $D \Sigma_+ \Sigma'$ with
$\Sigma' \in {\rm Aut}_{R_0(I)} (\Gamma_0)$.

Turning to (a), let $\Sigma$ be any morphism
in $R_0(I)$. By analogy with Coxeter groups, we define the {\em length} of
$\Sigma$, denoted $\ell(\Sigma)$, to be the minimum length of an
expression
$\Sigma=\Sigma_{i_1}\Sigma_{i_2}\cdots \Sigma_{i_r}$. We call an
expression
of minimal length for $\Sigma$ a {\em reduced expression} for $\Sigma$.
In the sequel, we will often write $i$ instead of $\Sigma_i$ to simplify
the
notation.

The following Lemma is an immediate consequence of the relations (R1) and
(R2):

\begin{lemma} \label{reduced}
An expression $\Sigma_{i_1}\Sigma_{i_2}\cdots \Sigma_{i_r}$ is reduced if
and only if between any two occurrences of the same $i$ there is an
occurrence
of $j$ for some $j$ linked to $i$.
\end{lemma}

We can strengthen this as follows:

\begin{lemma} \label{inbetween}
An expression $\Sigma_{i_1}\Sigma_{i_2}\cdots \Sigma_{i_r}$ for a morphism
in $R_0(I)$ is reduced if and only if between any two occurrences of
$i$ there is precisely one occurrence of $j$ for each $j$
linked to $i$.
\end{lemma}

\proof The ``if'' part follows from Lemma~\ref{reduced} so we only have to
show the ``only if'' statement. So suppose that
$\Sigma_{i_1}\Sigma_{i_2}\cdots \Sigma_{i_r}$ is a reduced expression.
We argue by induction on the distance between two consecutive occurrences
of
$i$. If the distance is $1$ then the corresponding interval is
$\Sigma_i\Sigma_j\Sigma_i$. Since this expression is reduced, $j$ is
linked
to $i$, and in order for the second $\Sigma_i$ to be applicable, $j$ must
be
the only element of $I$ that is linked to $i$, so we are done.

For the inductive step, we must first show that every $j$ linked to $i$
appears at most once between two consecutive occurrences of $i$. This is
true since, by induction, if we had two occurrences of some $j$ linked to
$i$, between them $i$ should occur, a contradiction. Then notice that all
$j$ linked to $i$ must appear between the two consecutive occurrences of
$i$,
because
we have to reverse all arrows at $i$.
\endproof

\begin{lemma} \label{extremal}
Suppose that $\Sigma_{i_1}\Sigma_{i_2}\cdots \Sigma_{i_p}$ is a reduced
expression for a morphism in $R_0(I)$ that starts at an arbitrary
orientation
of $I$ and ends up at $\Gamma_0$ or $\Gamma_0^{\rm op}$.
Then, for every $i$ that appears in the reduced expression, we have
precisely
one of the following:
\begin{itemize}
\item[(i)] no $j$ linked to $i$ appears before the first occurrence of
$i$,
or
\item[(ii)] every $j$ linked to $i$ appears exactly once before the first
occurrence of $i$.
\end{itemize}
\end{lemma}

\proof This follows easily from Lemma~\ref{inbetween}. Let
$\Sigma'=\Sigma_i \Sigma$; note that this composition is well-defined
since $\Sigma$ ends up at $\Gamma_0$ or $\Gamma_0^{\rm op}$,
from which every $\Sigma_i$ is
applicable. If $\Sigma_i\Sigma_{i_1}\Sigma_{i_2}\cdots \Sigma_{i_p}$ is
reduced then (ii) holds by Lemma~\ref{inbetween}, otherwise (i) holds by
Lemma~\ref{reduced}.
\endproof

\begin{lemma} \label{allappear}
Suppose that $\Sigma$ is a morphism of $R_0(I)$ which begins and ends at
the
same orientation $\Gamma$, and that $\Sigma$ is not the identity morphism.
Then every reduced expression for $\Sigma$ contains all $i\in I$.
\end{lemma}

\proof It is enough to show that if $i$ occurs in a reduced expression for
$\Sigma$ and $j$ is linked to $i$, then $j$ also occurs. Suppose this is
not
the case. Then $i$ must occur at least twice, because the arrow linking
$i$ and $j$ must be reversed at least twice by $\Sigma$.
So $j$ must occur by Lemma~\ref{inbetween}.
\endproof

\noindent
Now everything is ready for the proof of Theorem~\ref{alternatingaut}(a).
We need to show that every
morphism $\Sigma$ in $R_0(I)$ which begins and ends at the same quiver
$\Gamma$, where $\Gamma$ is $\Gamma_0$
or $\Gamma_0^{\rm op}$, has a reduced expression of form
$(\Sigma_+\Sigma_-)^k$ or $(\Sigma_-\Sigma_+)^k$ for some
$k\in\mathbb{Z}_{\geq 0}$.

We proceed by induction on $\ell(\Sigma)$. There is nothing to prove when
$\ell(\Sigma)=0$, so let us assume that $\Sigma$ is non-trivial, and that
the statement is known for morphisms of smaller length. Fix a reduced
expression for $\Sigma$. By Lemma~\ref{allappear}, all $i\in I$ occur.
By Lemma~\ref{extremal}, $I$ is the disjoint union of two subsets: those
satisfying~\ref{extremal}(i) and those satisfying~\ref{extremal}(ii).
By definition, if $i$ and $j$ are linked then they belong to different
subsets, since whichever of $i$ and $j$ appears first in
$\Sigma$ satisfies (i), while the other satisfies (ii).
It follows that one of these parts is $I^+$, and the other is $I^-$.
Suppose
that $I^+$ is the set of indices satisfying (i) (the argument if $I^-$ is
this set is similar). Commuting the first
occurrences of $\Sigma_i$ for each $i\in I^+$ to the beginning of our
reduced expression, we obtain a reduced expression for $\Sigma$ of the
form
$\Sigma_+\Sigma'$ with $\ell(\Sigma')=\ell(\Sigma)-|I^+|$, such that every
reduced expression for $\Sigma'$ starts with $\Sigma_j$ for some
$j\in I^-$. Note that $\Sigma'$ is a morphism from $\Gamma$ to
$\Gamma^{\rm op}$.

Next consider the morphism $\Sigma''=\Sigma'\Sigma_+$ from $\Gamma^{\rm
op}$
to $\Gamma^{\rm op}$. It is clear that
$\ell (\Sigma'') \leq \ell (\Sigma') + |I^+| = \ell (\Sigma)$.
If $\ell (\Sigma'') < \ell (\Sigma)$ then by induction
$\Sigma''$ has an expression $(\Sigma_+ \Sigma_-)^k$ or
$(\Sigma_- \Sigma_+)^k$, which implies the same claim for
$\Sigma = \Sigma_+\Sigma'=\Sigma_+ \Sigma'' \Sigma_+$.
or

So let us assume that $\ell (\Sigma'') = \ell (\Sigma)$ and take a
reduced expression for $\Sigma''$ obtained by taking the original reduced
expression for $\Sigma'$ (recall that it starts with some $j \in I^-$)
followed by any reduced expression for $\Sigma_+$.
Applying to this expression the same argument as before, we conclude that
by moving some $\Sigma_j$'s for $j \in I^-$ to the left we obtain a
reduced expression for $\Sigma'$ that starts with $\Sigma_-$.
This leads to a reduced expression for $\Sigma$ that starts with
$\Sigma_+ \Sigma_-$; by induction the rest of the reduced expression can
be written in the desired form, and Theorem~\ref{alternatingaut}(a)
follows.
\endproof

\section{An application to quiver representations}

As before, let $\Gamma$ be a quiver whose underlying graph $I$ is
a Dynkin diagram of the $ADE$-type.
We say that a set $C$ of isomorphism classes of
indecomposable representations of $\Gamma$ is
\emph{Ext-free} (also called a partial tilting set) if ${\rm Ext}^1 (M,N)
= 0$ for every $M, N \in C$.
For any nonnegative integer $k$, let $f_\Gamma^+ (k)$ denote the
number of these Ext-free sets of cardinality $k$.

\begin{proposition}
\label{pr:ext-free-enumeration}
For every $k$, the number $f_\Gamma^+ (k)$ depends only on the
underlying Dynkin graph, not on the choice of an orientation $\Gamma$.
\end{proposition}

\proof
Using Gabriel's theorem, we identify isoclasses of indecomposable
representations of $\Gamma$ with positive roots in $\Phi$.
Under this identification, the Ext-free sets become positive
$\Gamma$-compatible subsets of $\Phi_{\geq -1}$, in the
terminology of Definition~\ref{def:compatible}.

For every subset $J \subset I$, let $f_\Gamma (k,J)$
(resp. $f^+_\Gamma (k,J)$) denote the number of
\mbox{$\Gamma|_J$-compatible} subsets
of cardinality $k$ in $\Phi (J)_{\geq -1}$ (resp.
$\Phi (J)_{> 0}$); here $\Gamma|_J$ is the orientation induced by $\Gamma$
on the full subgraph $J$ of $I$.
Thus, we have in particular $f^+_\Gamma (k,I) = f^+_\Gamma (k)$.
In view of Proposition~\ref{redpos}, these numbers satisfy the following
relations:
$$f_\Gamma (k,J) = \sum_{K \subset J, |K| \leq k}
f^+_\Gamma (k-|K|,J-K) \ .$$
It follows that the numbers $f^+_\Gamma (k,J)$ can be expressed
as integer linear combinations of the numbers $f_\Gamma (l,K)$
by using the M\"obius inversion on the partially ordered set of
pairs $(k,J)$ with the partial order given by
$$(l,K) \leq (k,J) \Leftrightarrow K \subset J, k-l = |J - K| \ .$$
Since this poset structure is independent of $\Gamma$, it remains
to show that so are the numbers $f_\Gamma (k,I)$.
But the latter statement is an immediate consequence of
Corollary~\ref{cor:independence}.
\endproof

For every $\Gamma$, the corresponding Ext-free sets form a simplicial
complex,
which was studied by C.~Riedtmann and A.~Schofield~\cite{RS} and
L.~Unger~\cite{U} (see Remark~\ref{ungerremark}).
Proposition~\ref{pr:ext-free-enumeration} shows that the
$f$-vector of this complex depends only on the underlying Dynkin graph.
Despite this fact, the complexes corresponding to different orientations
of the same graph may be not isomorphic to each other.

\begin{figure}[ht]
\begin{center}
\setlength{\unitlength}{3pt}
\begin{picture}(100,60)(0,-10)
\thicklines
\put(0,0){\line(1,0){40}}
\put(0,0){\line(1,2){20}}
\put(20,40){\line(1,-2){20}}
\put(0,0){\line(3,2){30}}
\put(40,0){\line(-3,2){30}}
\put(10,20){\line(1,0){20}}

\multiput(0,0)(40,0){2}{\circle*{2}}
\multiput(10,20)(20,0){2}{\circle*{2}}
\put(20,40){\circle*{2}}
\put(20,13){\circle*{2}}

\put(0,-4){\makebox(0,0){$1$}}
\put(40,-4){\makebox(0,0){$3$}}
\put(6,20){\makebox(0,0){$12$}}
\put(34,20){\makebox(0,0){$23$}}
\put(20,9){\makebox(0,0){$123$}}
\put(20,44){\makebox(0,0){$2$}}
\put(20,-8){\makebox(0,0){$\Gamma_0$}}

\put(60,0){\line(1,0){40}}
\put(60,0){\line(1,2){20}}
\put(80,40){\line(1,-2){20}}
\put(60,0){\line(3,2){30}}
\put(100,0){\line(-3,2){30}}
\put(80,40){\line(0,-1){27}}

\multiput(60,0)(40,0){2}{\circle*{2}}
\multiput(70,20)(20,0){2}{\circle*{2}}
\put(80,40){\circle*{2}}
\put(80,13){\circle*{2}}

\put(60,-4){\makebox(0,0){$1$}}
\put(100,-4){\makebox(0,0){$3$}}
\put(66,20){\makebox(0,0){$12$}}
\put(94,20){\makebox(0,0){$23$}}
\put(80,9){\makebox(0,0){$123$}}
\put(80,44){\makebox(0,0){$2$}}
\put(80,-8){\makebox(0,0){$\Gamma_1$}}

\end{picture}
\end{center}
\caption{``Positive" simplicial complexes in type $A_3$}
\label{fig:a3}
\end{figure}
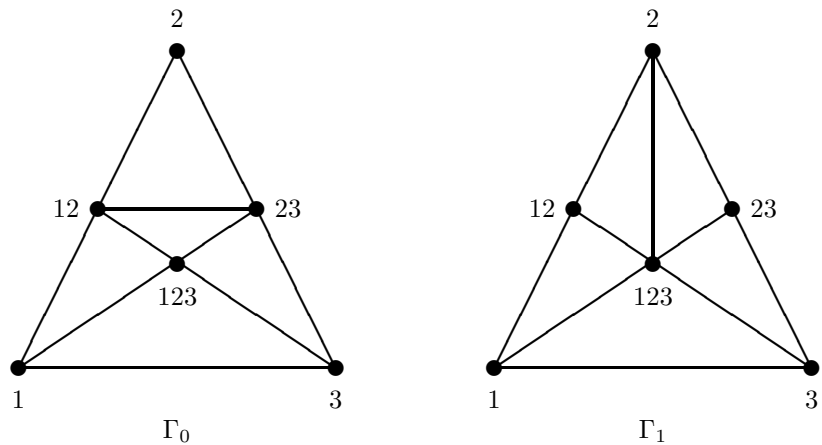

\begin{example}
{\rm Let $I$ be of type $A_3$, and let $\Gamma_0$ (resp. $\Gamma_1$)
be an alternating (resp. equioriented) orientation of $I$.
The corresponding simplicial complexes are shown in Figure~\ref{fig:a3},
and they are clearly non-isomorphic (in the figure, the vertices
labeled by $i, ij$, and $ijk$ represent roots $\alpha_i$,
$\alpha_i + \alpha_j$, and $\alpha_i + \alpha_j + \alpha_k$,
respectively).
}
\end{example}

A recursive way to compute the numbers $f_\Gamma (k,I)$ is given
by \cite[Proposition~3.7]{FZ}.
The numbers $f^+_\Gamma (k,I)$ can then be found by the M\"obius
inversion described in the proof of
Proposition~\ref{pr:ext-free-enumeration}.
In particular, when $\Phi$ is irreducible, the total number
$f^+_\Gamma (|I|,I)$ of positive
$\Gamma$-clusters is given by \cite[(3.8)]{FZ}:
$$f^+_\Gamma (|I|,I) = \prod_{i \in I} \frac{e_i + h - 1}{e_i + 1} \ ,$$
where the $e_i$ are the exponents of $\Phi$, and $h$ is the Coxeter
number.
It would be interesting to explain this formula using the theory
of quiver representations.

\end{document}